%% file: JewishProblems.tex
\documentclass{article}
\usepackage[pdftex]{graphicx}
\usepackage{amsmath}
\usepackage{amsfonts}
\usepackage{amsthm}
\usepackage{verbatim}
\usepackage{color}
\usepackage{ifpdf}
\usepackage{url}

\ifpdf
\else
\fi

\newcounter{problem}
\newcommand{\problem}{\refstepcounter{problem}\subsection*{Problem~\theproblem}}

\newcommand{\pt}[1]{{{#1}})}

\newcounter{cases}
\newcommand{\case}[1]{\mbox{}\\\noindent{\bf Case \thecases}, {{#1}}:\\%
\addtocounter{cases}{1}}
\newcommand{\stopcases}{\setcounter{cases}{1}}

\newcommand{\ten}{10^\circ}

\begin{document}

\title{Jewish Problems}

\author{Tanya Khovanova\\Massachusetts Institute of Technology\\Cambridge, MA, USA \and Alexey Radul\\Hamilton Institute at NUIM\\Maynooth, Co.~Kildare, Ireland}

\maketitle

\begin{abstract}
This is a special collection of problems that were given to select applicants during oral entrance exams to the math department of Moscow State University. These problems were designed to prevent Jewish people and other undesirables from getting a passing grade. Among problems that were used by the department to blackball unwanted candidate students, these problems are distinguished by having a simple solution that is difficult to find. Using problems with a simple solution protected the administration from extra complaints and appeals. This collection therefore has mathematical as well as historical value.
\end{abstract}

\section{A personal story of Tanya Khovanova}

In the summer of 1975, while I was in a Soviet math camp preparing to compete in the International Math Olympiad on behalf of the Soviet Union, my fellow team members and I were approached for help by Valera Senderov, a math teacher in one of Moscow's best special math schools.

The Mathematics Department of Moscow State University, the most prestigious mathematics school in Russia, was at that time actively trying to keep Jewish students (and other ``undesirables'') from enrolling in the department. One of the methods they used for doing this was to give the unwanted students a different set of problems on their oral exam. I was told that these problems were carefully designed to have elementary solutions (so that the Department could avoid scandals) that were nearly impossible to find. Any student who failed to answer could easily be rejected, so this system was an effective method of controlling admissions. These kinds of math problems were informally referred to as ``Jewish'' problems or ``coffins''. ``Coffins'' is the literal translation from Russian; they have also been called ``killer'' problems in English.

These problems and their solutions were, of course, kept secret, but Valera Senderov and his friends had managed to collect a list. In 1975, they approached us to solve these problems, so that they could train the Jewish and other students in these mathematical ideas. Our team of the best eight Soviet students, during the month we had the problems, solved only half of them. True, that we had other priorities, but this fact speaks to the difficulty of these problems.

Being young and impressionable, I was shaken by this whole situation.  I had had no idea that such blatant discrimination had been going on.  In addition to trying to solve them at the time, I kept these problems as my most valuable possession---I still have that teal notebook.

Later, I emigrated to the United States.  When I started my own web
page, one of the first things I did was to post some of
the problems. People sent me more problems, and solutions to the ones I
had.  It turned out that not all of the coffins even had elementary
solutions: some were intentionally ambiguous questions, some were just
plain hard, some had impossible premises. This article is a selection
from my collection; we picked out some choice problems that do contain
interesting tricks or ideas.

\hfill Tanya Khovanova

\section{Introduction}

Discrimination against Jewish people at the entrance exams to the most prestigious universities in the USSR is a documented fact, see A.~Shen~\cite{Shen} and A.~Vershik~\cite{Vershik}. Alexander Shen in his article published 25 problems that were given to Jewish applicants. Later Ilan Vardi wrote solutions to all those problems, see I.~Vardi's homepage~\cite{Vardi}. These and other articles were published later in a book, ``You Failed Your Math Test, Comrade Einstein: Adventures and Misadventures of Young Mathematicians, Or Test Your Skills in Almost Recreational Mathematics.''~\cite{Shifman}.

In the present collection we have only one problem that overlaps with the problems that Vardi solved, and we give an easier and more surprising solution.

Now, after thirty years, these problems seem easier. Mostly, this is because the ideas of how to solve these problems have spread and are now a part of the standard set of ideas. Thirty years ago these problems were harder to solve and, in addition, the students were given these problems one after another until they failed one of them, at which point they were given a failing mark.

To give readers an opportunity to try and solve these problems we separated problems, key ideas, and solutions. Section~\ref{sec:problems} contains 21 problems. In the subsequent Section~\ref{sec:ideas} we give hints or main ideas for the solutions. Finally, Section~\ref{sec:solutions} contains solutions.

\section{Problems}\label{sec:problems}

\problem \label{pr:1} 

Solve the following inequality for positive $x$:
\begin{equation}
 x(8 \sqrt{1-x} + \sqrt{1+x}) \leq 11 \sqrt{1+x} - 16 \sqrt{1-x}. 
   \label{pr1:eq} 
\end{equation}

\problem \label{pr:7} 

Find all functions $F(x): \mathbb{R} \rightarrow \mathbb{R}$ having the property that for any $x_1$ and $x_2$
the following inequality holds:
\begin{equation} F(x_1) - F(x_2) \leq (x_1 - x_2)^2. \label{pr7:prob} 
\end{equation}

\problem \label{pr:10} 

Given a triangle $ABC$, construct, using straightedge and compass, a point $K$ 
on $AB$ and a point $M$ on $BC$, such that 
\begin{equation} AK = KM = MC. \label{pr10:prob} \end{equation}

\problem \label{pr:12} 

Solve the following equation for real $y$:
\begin{equation} 2\sqrt[3]{2y-1} = y^3 + 1. \label{pr12:prob} \end{equation}

\problem \label{pr:17} 

Solve the equation
\begin{equation} \sin^7 x + \frac{1}{\sin^3 x} = \cos^7 x + \frac{1}{\cos^3 x}.
\label{pr17:prob} \end{equation}

\problem \label{pr:19} 

You are given a point $M$ and an angle $C$ in the plane.  Using a ruler and a compass, draw a line through the point $M$ which cuts off the angle $C$ a 
triangle of

\pt{a} a given perimeter $p$

\pt{b} minimum perimeter.

\problem \label{pr:22} 

There is a circle in the plane with a drawn diameter.  Given a point, draw the
perpendicular from the point to the diameter using only a straightedge.  Assume the 
point is neither on the circle nor on the diameter line.

\problem \label{pr:26} 
Given an equilateral triangle $ABC$ and a point $O$ inside it, with $\angle BOC = x$ and $\angle AOC = y$, find, in terms of $x$ and $y$, 
the angles of the triangle with side lengths equal to $AO$, $BO$, and $CO$.

\problem \label{pr:30} 
Two intersecting lines are given on a plane.  Find the locus of points $A$
such that the sum of the distances from $A$ to each line is equal to a 
given number.

\problem \label{pr:31} 
A quadrilateral is given in space, such that its edges are tangent to a 
sphere.  Prove that all the points of tangency lie in one plane.

\problem \label{pr:38} 
Prove that $\sin 10^\circ$ is irrational.

\problem \label{pr:42} 
Can you put six points on the plane, so that the distance between any 
two of them is an integer, and no three are collinear?

\problem \label{pr:45} 
Is it possible to put an equilateral triangle onto a square grid so that all the vertices are in corners?

\problem \label{pr:48} 
\begin{figure}[htbp]
\begin{center}
\ifpdf
  \input{pr48.pdftex_t}
\else
  \input{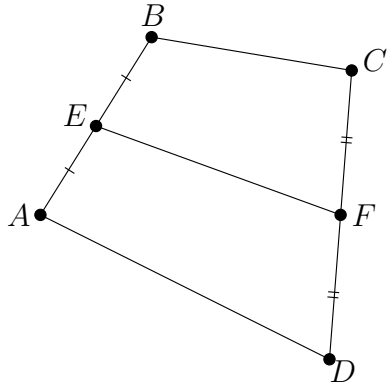}
\fi
\end{center}
\caption{Problem~\ref{pr:48}: Reconstruct this picture with ruler and
  compass, given only the lengths $|AB|, |BC|, |CD|, |DA|, |EF|$.}
\label{pr48:fig}
\end{figure}

Reconstruct the quadrilateral $ABCD$ given the lengths of its sides (in order) and the length of the midline between the first and third sides (namely, all the segments drawn in Figure~\ref{pr48:fig}).

\problem \label{pr:49} 
Find the quadrilateral with the largest area given the lengths of its four sides (in order).

\problem \label{pr:50} 
You are given two parallel segments. Using a straightedge, divide one of them into six equal parts.

\problem \label{pr:52} 
What's larger, $\log_2 3$ or $\log_3 5$?

\problem \label{pr:61} 
How many digits does the number $125^{100}$ have?

\problem \label{pr:65} 
Get rid of the radicals in the denominator of:
\[
\frac{1}{\sqrt[3]{a}+\sqrt[3]{b}+\sqrt[3]{c}}.
\]

\problem \label{pr:68} 
Construct (with ruler and compass) a square given one point from each side.

\begin{figure}[htbp]
\begin{center}
\ifpdf
  \input{pr71.pdftex_t}
\else
  \input{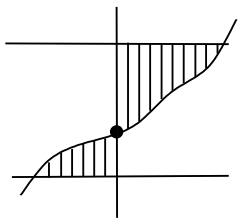}
\fi
\end{center}
\caption{Problem~\ref{pr:71}: Given the graph, find the point that
  minimizes the depicted area.}
\label{pr71:fig}
\end{figure}

\problem \label{pr:71} 
The graph of a monotonically increasing function is cut off with two horizontal lines. Find a point on the graph between intersections such that the sum of the two areas bounded by the lines, the graph and the vertical line through this point is minimum. See Figure~\ref{pr71:fig}.

\section{Ideas}\label{sec:ideas}

\subsection*{Idea for problem~\ref{pr:1}}
Substitute
\begin{equation} y = \frac{\sqrt{1-x}}{\sqrt{1+x}}. \label{pr1:sub} 
\end{equation}

\subsection*{Idea for problem~\ref{pr:7}}
Use derivatives.

\subsection*{Idea for problem~\ref{pr:10}}

Construct a triangle $A'B'C$ similar to $ABC$ around points $K'$ and $M'$ 
so that  the desired property holds for it.

\subsection*{Idea for problem~\ref{pr:12}}
Define
\[ x = (y^3 + 1)/2. \]

\subsection*{Idea for problem~\ref{pr:17}}
$\cos^{2k+1}x - \sin^{2k+1}x = (\cos x - \sin x)P(\sin x \cos x)$, where $P$ is a polynomial.

\subsection*{Idea for problem~\ref{pr:19}}
Inscribe a circle into the angle $C$.

\subsection*{Idea for problem~\ref{pr:22}}
\begin{figure}[htbp]
\begin{center}
\ifpdf
  \input{pr22.pdftex_t}
\else
  \input{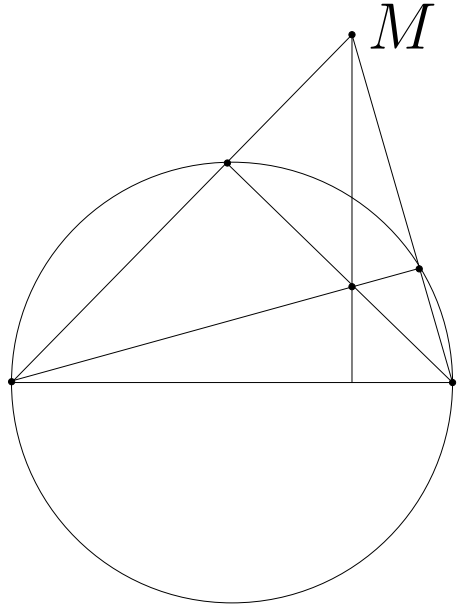}
\fi
\end{center}
\caption{Idea for Problem~\ref{pr:22}.}
\label{pr22:fig}
\end{figure}

See Figure~\ref{pr22:fig}.

\subsection*{Idea for problem~\ref{pr:26}}
\begin{figure}[htbp]
\begin{center}
\ifpdf
  \input{pr26.pdftex_t}
\else
  \input{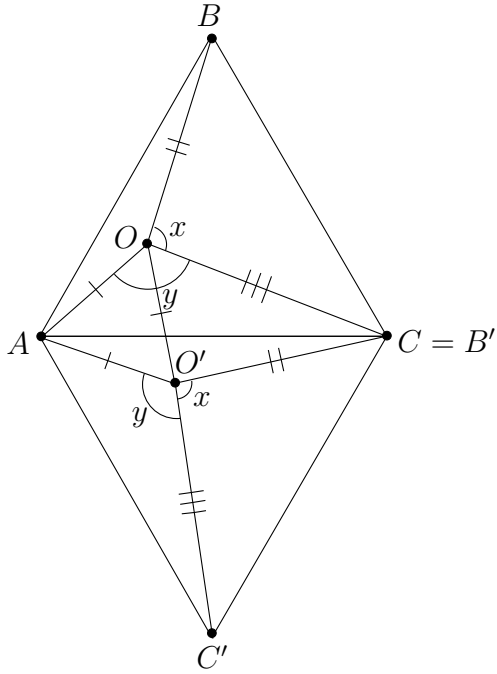}
\fi
\end{center}
\caption{Rotation for Problem~\ref{pr:26}.}
\label{pr26:fig}
\end{figure}

Rotate the picture by $60^\circ$ 
about $A$.  Let the image of $O$ be $O'$.
Then the triangle $COO'$ (or the triangle $BOO'$ if we rotate in the other direction) will be 
the desired one. See Figure~\ref{pr26:fig}.

\subsection*{Idea for problem~\ref{pr:30}}
In an isosceles triangle the sum of the distances from any point on the base to the two other sides is fixed.

\subsection*{Idea for problem~\ref{pr:31}}
Use gravity.

\subsection*{Idea for problem~\ref{pr:38}}
Express $\sin 10^\circ$ through $\sin 30^\circ$.

\subsection*{Idea for problem~\ref{pr:42}}
A Pythagorean triangle generates three such points.  Various reflections of it can increase the number of points.

\subsection*{Idea for problem~\ref{pr:45}}
Use parity considerations or the fact that $\tan 60^\circ$ is irrational.

\subsection*{Idea for problem~\ref{pr:48}}
\begin{figure}[htbp]
\begin{center}
\ifpdf
  \input{pr48sol1.pdftex_t}
\else
  \input{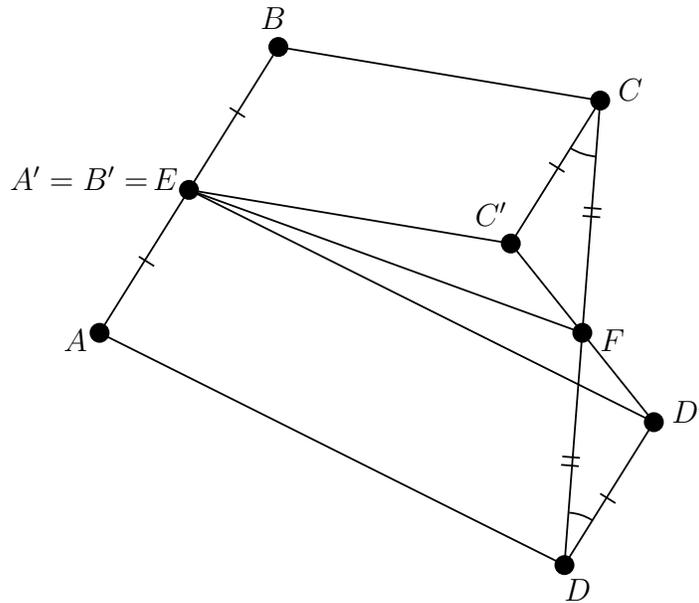}
\fi
\end{center}
\caption{Parallel translations for Problem~\ref{pr:48}.}
\label{pr48sol1:fig}
\end{figure}

Consider the parallel translation of the sides $AD$ and $BC$ so that
$A'$ and $B'$ are at $E$, shown in Figure~\ref{pr48sol1:fig}.  The midline will become the median of the new triangle $EC'D'$.

\subsection*{Idea for problem~\ref{pr:49}}
It's the cyclic quadrilateral.

\subsection*{Idea for problem~\ref{pr:50}}
The construction can be separated into two tasks: divide a segment into two equal parts, and divide a segment into six equal parts if six equal segments on a parallel line are given.

\subsection*{Idea for problem~\ref{pr:52}}
Compare both values to 3/2.

\subsection*{Idea for problem~\ref{pr:61}}
Use the fact that $2^{10} = 1024$, which is close to $10^3$.

\subsection*{Idea for problem~\ref{pr:65}}
Use complex numbers.

\subsection*{Idea for problem~\ref{pr:68}}
Cut two parallel lines at a given distance by another line at a given angle. The length of the cut out segment is fixed.

\subsection*{Idea for problem~\ref{pr:71}}
The point on the middle line in between the two lines should work.

\section{Solutions}\label{sec:solutions}

\subsection*{Solution to problem~\ref{pr:1}}
First, observe that for $x > 1$, the terms in (\ref{pr1:eq}) become undefined.
Next, define $y$ per the idea (\ref{pr1:sub}):
\[ y = \frac{\sqrt{1-x}}{\sqrt{1+x}}.
\]
Observe that for permissible
values of $x$, we have
\[ 0 \leq y \leq 1, \]
and $y$ is monotonically decreasing in $x$.
Observe also that
\[ x = \frac{1-y^2}{1+y^2}. \label{pr1:subinv} \]

Given the above preliminaries, our inequality transforms as follows:
\begin{eqnarray*}
 x(8 \sqrt{1-x} + \sqrt{1+x}) & \leq & 11 \sqrt{1+x} - 16 \sqrt{1-x} \\
 x \left( 8 \frac{\sqrt{1-x}}{\sqrt{1+x}} + 1 \right) & \leq & 
    11 - 16 \frac{\sqrt{1-x}}{\sqrt{1+x}} \\
 \frac{1-y^2}{1+y^2}(8y + 1) & \leq & 11 - 16y \\
 (1 - y^2)(8y + 1) & \leq & (1 + y^2)(11 - 16y) \\
 -8y^3 - y^2 + 8y + 1 & \leq & -16y^3 + 11y^2 -16y + 11 \\
 -8y^3 + 12y^2 - 24y + 10 & \geq & 0 \\
 (2y-1)(-4y^2 + 4y - 10) & \geq & 0.
\end{eqnarray*}
Now, $-4y^2 + 4y - 10$ is always negative, so our monstrous inequality 
reduces, in the end, to the humble
\[ 2y - 1 \leq 0. \]
Invoking the monotonicity of $y$ in $x$ and the fact that
\[ \frac{1 - (1/2)^2}{1 + (1/2)^2} = \frac{3}{5}, \]
our final answer is
\[ \frac{3}{5} \leq x \leq 1. \]

\subsection*{Solution to problem~\ref{pr:7}}
Inequality (\ref{pr7:prob}) implies
\[ \frac{F(x_1) - F(x_2)}{|x_1 - x_2|} \leq |x_1 - x_2|, \]
so the derivative of $F$ at any point $x_2$ exists and is equal to zero.
Therefore, by the fundamental theorem of calculus, the constant functions are exactly the functions with the 
desired property.

\subsection*{Solution to problem~\ref{pr:10}}
\begin{figure}[htbp]
\begin{center}
\ifpdf
  \input{pr10.pdftex_t}
\else
  \input{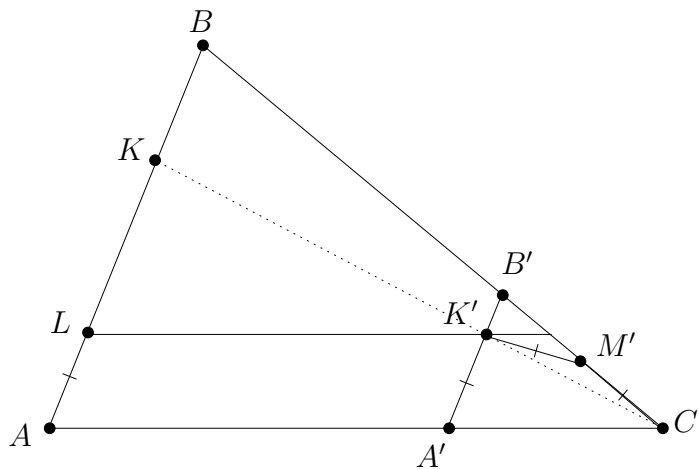}
\fi
\end{center}
\caption{Construction for Problem~\ref{pr:10}, in the order $M', L, K', K$.}
\label{pr10:fig}
\end{figure}

Pick a point $M'$ anywhere on $BC$.  Construct $L$ on $AB$ such that 
\[ AL = CM'. \]
Now we want to parallel translate $AL$ so that its image will satisfy
(\ref{pr10:prob}) with $CM'$. To do this, construct through $L$ a line
parallel to $AC$, and on it construct the point $K'$ such that
\[ K'M' = CM'. \]
Construct through $K'$ a line parallel to $AB$, and let it intersect
the angle $ACB$ at points $A'$ and $B'$.  Then the triangle $A'B'C$
and the points $K'$ and $M'$ satisfy (\ref{pr10:prob}), and $A'B'C$ is
similar to $ABC$.  Perform a homothety (also called a dilation) about
$C$ and you're done.  Figure~\ref{pr10:fig} illustrates.

\subsection*{Solution to problem~\ref{pr:12}}

Define
\[ x = (y^3 + 1)/2. \]
Then (\ref{pr12:prob}) becomes
\[ y = (x^3 + 1)/2. \]
Since $x$ and $y$ are given by the same function $g$ of each other, and 
$g$ is monotonically increasing, we can conclude
\[ x = y. \]
The equation is now a standard high-school cubic
\[ y^3 -2y +1 = 0, \]
solvable by the method of guessing a small integer root and factoring. Indeed,
\[ y = 1 \]
is a root, and the factor $y^2+y-1$ yields two more roots
\[ y = \frac{-1 \pm \sqrt{5}}{2}. \]

\subsection*{Solution to problem~\ref{pr:17}}
First, rearrange (\ref{pr17:prob}) into
\begin{eqnarray*} 
 \frac{1}{\sin^3 x} - \frac{1}{\cos^3 x} & = & \cos^7 x - \sin^7 x \\
 \frac{\cos^3 x - \sin^3 x}{\cos^3 x \sin^3 x} & = &
 \cos^7 x - \sin^7 x. 
\end{eqnarray*}
Now there are two cases:
\setcounter{cases}{1}
\case{$\cos x - \sin x = 0$}
Then
\begin{eqnarray*} x = \frac{\pi}{4} & \textrm{ or } & x = \pi + \frac{\pi}{4}.
\end{eqnarray*}
\case{otherwise}
We can cancel $\cos x - \sin x$, and get
\begin{align}
 \lefteqn{\frac{\cos^2 x + \cos x \sin x + \sin^2 x}{\cos^3 x \sin^3 x}  = } \notag \\
  & &  \cos^6 x + \cos^5 x \sin x + \cos^4 x \sin^2 x  + \cos^3 x \sin^3 x + \notag \\
& & \cos^2 x \sin^4 x + \cos x \sin^5 x + \sin^6 x. \label{sol17:big}
\end{align}
Substituting 
\[ t = \cos x \sin x, \]
we find that
\begin{eqnarray*} 
 \cos^4 x + \sin^4 x & = & (\cos^2 x + \sin^2 x)^2 - 2\cos^2 x \sin^2 x \\
 & = & 1 - 2t^2 \end{eqnarray*}
and that
\begin{eqnarray*}
 \cos^6 x + \sin^6 x & = & (\cos^2 x + \sin^2 x)(\cos^4 x - \cos^2 x \sin^2 x 
    + \sin^4 x) \\
 & = & 1 - 3t^2 \end{eqnarray*}
so (\ref{sol17:big}) reduces to
\begin{eqnarray*}
 \frac{1 + t}{t^3} & = & 1 - 3t^2 + t(1 - 2t^2) + t^2 + t^3 = 1 + t - 2t^2 
   - t^3 \\
 0 & = & -t^6 - 2t^5 + t^4 + t^3 - t - 1. \end{eqnarray*}
But
\[ |t| = \left|\frac{\sin 2x}{2}\right| \leq \frac{1}{2}, \]
so
\[ |-t^6 - 2t^5 + t^4 + t^3 - t| \leq
 \frac{1}{64} + \frac{1}{16} + \frac{1}{16} + \frac{1}{8} + \frac{1}{2} < 1.
\]
Therefore, there are no solutions in this case.
\stopcases

\subsection*{Solution to problem~\ref{pr:19}}
\begin{figure}[htbp]
\begin{center}
\ifpdf
  \input{pr19.pdftex_t}
\else
  \input{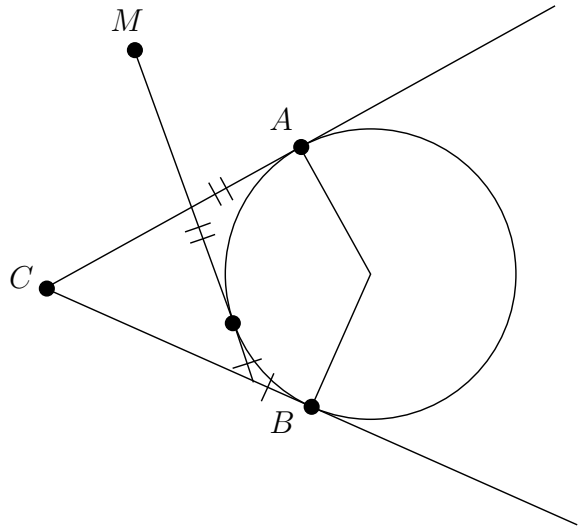}
\fi
\end{center}
\caption{Construction for Problem~\ref{pr:19}a.}
\label{pr19:fig1}
\end{figure}

\pt{a} The point $M$ must be outside the angle opposite to $C$. Construct points $A$ and $B$ on the angle such that 
\[ \frac{p}{2} = CA = CB. \]
Then inscribe a circle into the angle tangent at those two points, and draw 
a tangent from $M$ to that circle. The resulting triangle has the desired perimeter. This construction is illustrated in 
Figure~\ref{pr19:fig1}.

\pt{b}
If $M$ is outside the angle, the minimum perimeter is $0$.  If it is inside,
then draw a circle through $M$ away from $C$ that is tangent to the angle, and take a 
tangent to it as in part a.  It is not hard to see from Figure~\ref{pr19:fig1}
that the resulting triangle has minimum perimeter.

\subsection*{Solution to problem~\ref{pr:22}}
The end points of the diameter and $M$ form a triangle in which we need to draw a height, see Figure~\ref{pr22:fig}. The height passes through the orthocenter, which is easy to construct as the bases of the other two heights lie on the circle. If $M$ is inside the circle, the construction is the same.

\subsection*{Solution to problem~\ref{pr:26}}
The situation is represented in Figure~\ref{pr26:fig}.  The points $B$, $C$, 
and $O$ go to the points $B'$, $C'$, and $O'$, respectively,
under rotation about $A$ by $60^\circ$.  The marked equalities follow by 
the properties of rotation.  The desired triangle is then
$COO'$.  The angles in question are
\begin{eqnarray*} 
 O' O C & = & y - 60^\circ, \\
 OO' C & = & 300^\circ - x - y, \\
 OCO' & = & 180^\circ - (y - 60^\circ) - (300^\circ - x - y) \\
 & = & x - 60^\circ.
\end{eqnarray*}

\subsection*{Solution to problem~\ref{pr:30}}
\begin{figure}[htbp]
\begin{center}
\ifpdf
  \input{pr30i.pdftex_t}
\else
  \input{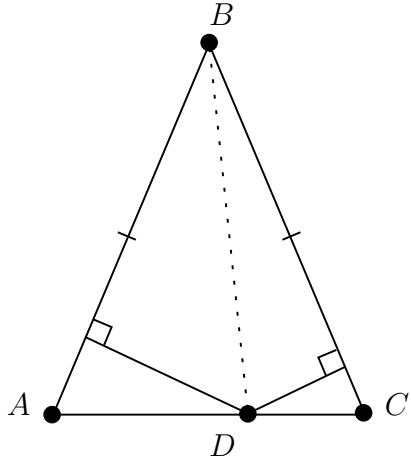}
\fi
\end{center}
\caption{One quadrant for Problem~\ref{pr:30}.  $|BC| = |AB|$, so the
  sum of the heights from $D$ is independent of the position of $D$ on
  $AC$.}
\label{pr30i:fig}
\end{figure}

In Figure~\ref{pr30i:fig}, the sum of the areas of the two triangles
$ABD$ and $CBD$ equals the area of the triangle $ABC$.  Canceling the
length of the side $BC (= AB)$, the sum of the heights from point
$D$ equals the height from point $A$ (and the height from point $C$).

\begin{figure}[htbp]
\begin{center}
\ifpdf
  \input{pr30.pdftex_t}
\else
  \input{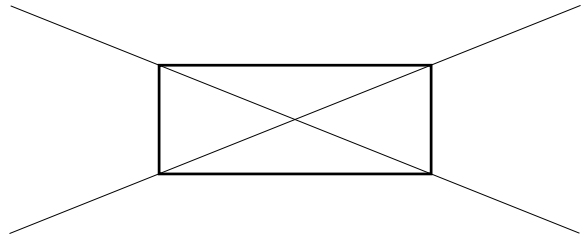}
\fi
\end{center}
\caption{Answer to Problem~\ref{pr:30}: locus of points with fixed
  total distance from two given lines.}
\label{pr30:fig1}
\end{figure}

Therefore, the locus is a rectangle whose vertices are the four points
that are each on one of the lines and at the desired distance from the
other: Figure~\ref{pr30:fig1}.

\subsection*{Solution to problem~\ref{pr:31}}
Place masses at the vertices of the quadrilateral such that the center 
of mass of each edge falls onto a point of tangency. One way to 
do this is to place, at each vertex, a mass equal to 1 over the distance 
from that the vertex to each point of tangency of the edges that vertex is on (the two distances must be equal). Then the center of mass of the whole arrangement will have to be 
on each line connecting pairs of opposite points of tangency.  Therefore,
these lines intersect, and all four points of tangency are in one plane,
as desired.

\subsection*{Solution to problem~\ref{pr:38}}
Observe that
\begin{eqnarray}
1/2 \ = \ \sin 30^\circ & = & 3\sin \ten - 4 \sin^3 \ten \label{pr38:sol1} \\
 0 & = & 8\sin^3 \ten - 6\sin \ten + 1, \label{pr38:sol2}
\end{eqnarray}
where (\ref{pr38:sol1}) holds by repeated application of the sine and cosine 
angle sum formulae.  Let
\[ x = 2\sin \ten. \]
Then (\ref{pr38:sol2}) reduces to 
\[ x^3 - 3x + 1 = 0. \]
All rational roots of this must be integers and divide the constant 
term.  Since $\pm 1$ do not work, all roots must be irrational.

\subsection*{Solution to problem~\ref{pr:42}}
If we
take a right triangle $ABC$ with integer sides, the area will be
rational.  Therefore, the height onto the hypotenuse $AC$ will be rational.  
This means that if we reflect this triangle about the hypotenuse, we will
have four points in the plane all at rational distances from each other.

Now, let the base of that height be at point $D$.  Then, since the
triangle $ABC$ is similar to the triangle $BDC$, the distances $AD$
and $CD$ will both be rational.  Therefore, reflecting the whole thing
about the perpendicular bisector of $AC$ yields two more points, so
that all six are at rational distances from each other.  All that's
left now is to expand by a multiplicative factor that clears all the
denominators.

Start with the 3, 4, 5 Pythagorean triangle, lay the hypotenuse along
the $x$ axis, with the center at the origin, perform the above
construction, and scale by 10 to clear denominators.  Six such points
are:
\[ (\pm 25, 0),\ (\pm 7, \pm 24). \]

\begin{figure}[htbp]
\begin{center}
\ifpdf
  \input{pr42a.pdftex_t}
\else
  \input{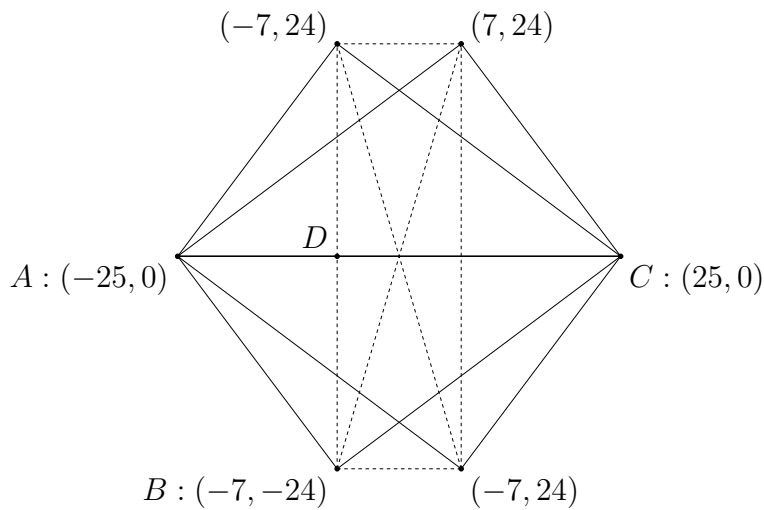}
\fi
\end{center}
\caption{Answer to Problem~\ref{pr:42}.  The six needed points have the
  given coordinates, and the points mentioned in the text are labeled.
  The solid lines have integral lengths by construction, being images
  under reflection of the initial Pythagorean triangle.  The dashed
  lines have lengths that must be proven or verified to be integral.}
\label{pr42a:fig}
\end{figure}

\subsection*{Solution to problem~\ref{pr:45}}
Suppose one of the vertices is at the origin. Suppose the coordinates of the two other vertices are $(a,b)$ and $(c,d)$. If all the numbers are divisible by two, we can reduce the triangle by half. Thus, we can assume that one of the numbers is not divisible by 2. Suppose $a$ is odd. 

If $b$ is odd, then $a^2 + b^2$ is of the form $4k+2$. As $a^2+b^2 = c^2 + d^2$ we must conclude that both $c$ and $d$ are odd. But the square of the length of the third side is $(a-c)^2 + (b-d)^2$, which is divisible by 4. So the triangle can't be equilateral.

If $b$ is even, then $a^2 + b^2$ is of the form $4k+1$. As $a^2+b^2 = c^2 + d^2$ we must conclude that $c$ and $d$ have different parity. Hence, $a-c$ and $b-d$ have the same parity and, correspondingly, the square of the length of the third side is again even. Contradiction.

For another solution, observe that the angle any segment between two
grid vertices makes with either axis has to have a rational
tangent. In addition, as
\[ \tan(\alpha + \beta) = \frac{\tan \alpha + \tan \beta}{1 - \tan
  \alpha \tan \beta}, \]
the tangent of the sum (or difference) of two angles with rational
tangents is rational.  Therefore any non-right angle given by
grid points has rational tangent.  Since
$\tan 60^\circ = \sqrt{3}$ is irrational, the equilateral triangle
cannot have all three vertices on a square grid.

\subsection*{Solution to problem~\ref{pr:48}}
Consider the parallel translation of the sides $AD$ and $BC$ so that
$A'$ and $B'$ are at $E$, shown in Figure~\ref{pr48sol1:fig}.  Then the
triangles $CC'F$ and $DD'F$ are equal, because they have two equal
sides with an equal angle between them.  Therefore, $F$ is the median
of the triangle $EC'D'$.  We can reconstruct a triangle from the
lengths of its two sides and the median between them,\footnote{Because
  we can reconstruct the parallelogram whose sides are the given sides
  and one of whose diagonals has twice the length of the given
  median.} which gives us the angle the original sides $AD$ and $BC$
made with each other.

\begin{figure}[htbp]
\begin{center}
\ifpdf
  \input{pr48sol2.pdftex_t}
\else
  \input{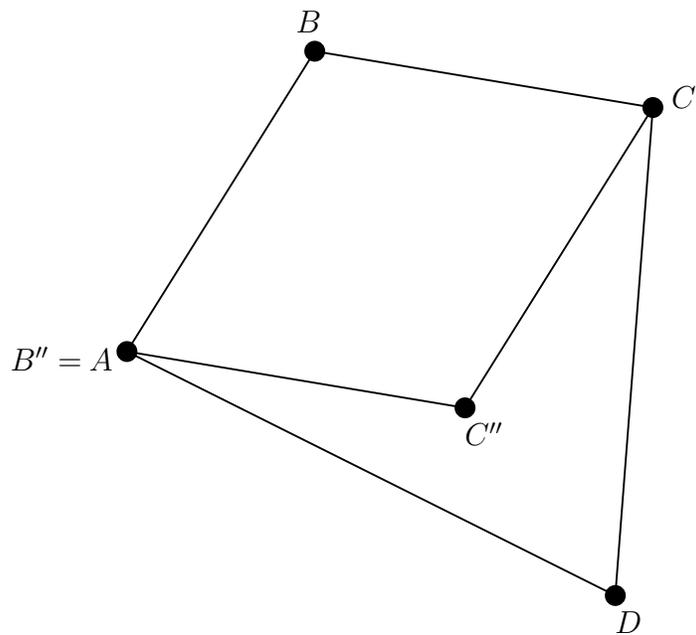}
\fi
\end{center}
\caption{Completion of Problem~\ref{pr:48}.}
\label{pr48sol2:fig}
\end{figure}

The reconstruction of that triangle is equal to the triangle $AC''D$
in Figure~\ref{pr48sol2:fig}, where $AC''$ is the parallel translation
of $BC$ so that $B''$ is at $A$.  But we are given the distance between
$C$ and $D$ and between $C$ and $C''$, so we can find $C$ by
constructing the intersection of two circles.  Constructing $B$ is
easy.

\subsection*{Solution to problem~\ref{pr:49}}
\begin{figure}[htbp]
\begin{center}
\ifpdf
  \input{pr49.pdftex_t}
\else
  \input{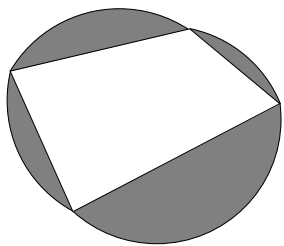}
\fi
\end{center}
\caption{A quadrilateral with attached parts, per the argument for Problem~\ref{pr:49}.}
\label{pr49:fig}
\end{figure}

Take the cyclic quadrilateral with the given sides and inscribe it into a circle. Consider the parts of the circle between the circumference and the quadrilateral as firmly attached to the sides. Consider a different quadrilateral with the given sides and the same parts attached to it shown in gray in Figure~\ref{pr49:fig}. The new figure will have the same perimeter as the perimeter of the circle, hence necessarily lesser area. Since the attached parts have fixed area, the quadrilateral must have lesser area as well.

\subsection*{Solution to problem~\ref{pr:50}}
Given a segment and a parallel line we can always divide the segment into two equal parts. Take two points on the parallel line. Together with the ends of the segment the points form a trapezoid. Continue the sides to form a triangle. The segment passing through the third vertex of the triangle and the intersection of the diagonals of the trapezoid divides both parallel segments into two equal parts.

\begin{figure}[htbp]
\begin{center}
\ifpdf
  \input{pr50.pdftex_t}
\else
  \input{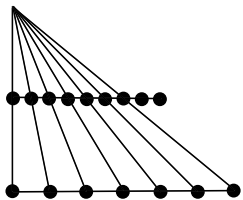}
\fi
\end{center}
\caption{The homothety for Problem~\ref{pr:50}.}
\label{pr50:fig}
\end{figure}

Using the division method above, divide one of the segments into eight equal parts. Pick six of these parts consecutively. Then perform a homothety mapping their union onto the other segment. The center of the homothety is the intersection of the other sides of the trapezoid formed by the six segments and the target segment. See Figure~\ref{pr50:fig}.

\subsection*{Solution to problem~\ref{pr:52}}
We will compare $\log_2 3$ and $\log_3 5$ to $3/2$, both by the method
of applying the same monotonic functions to both sides until we get
something known.

First, compare $\log_2 3$ and $3/2$. The functions $2^x$ and $x^2$ are
monotonic (at least for positive $x$), so we can use the composition of them to
remove ugly pieces from our numbers:

\[ \left(2^{\log_2 3}\right)^2 = 3^2 = 9 > 8 = 2^3 = \left(2^{3/2}\right)^2. \]

Second, compare $\log_3 5$ and $3/2$. The functions $3^x$ and $x^2$
are monotonic:

\[ \left(3^{\log_3 5}\right)^2 = 5^2 = 25 < 27 = 3^3 = \left(3^{3/2}\right)^2. \]

Hence $\log_2 3 > \log_3 5$.

\subsection*{Solution to problem~\ref{pr:61}}
Clearly, $125^{100} = 10^{300}/2^{300}$. Also, $2^{300} = 1024^{30} = 10^{90} \times 1.024^{30}$. So
\[ 125^{100} = 10^{210}/1.024^{30}. \]
Let's estimate $1.024^{30}$. For any number slightly bigger than 1 we can use the binomial formula: $(1+x)^{30} = 1 + 30x + 435 x^2 + \binom{30}{3} x^3 + \ldots$. In our case $x < 1/40$. We can estimate each term after 1 as 0.75, 0.27, and 0.06. We see that they are getting smaller very fast. So we know that $1 < 1.024^{30} < 10$, which gives us 210 digits for $125^{100}$.

\subsection*{Solution to problem~\ref{pr:65}}
Let $x = \sqrt[3]{a}$, $y =\sqrt[3]{b}$, and $z=\sqrt[3]{c}$.
Consider three numbers: 1, $w$, $w^2$, where $w$ is a primitive cube root of 1. Consider the nine terms of the form  $x + w^i y + w^j z$, where $i$ and $j$ are 0, 1, or 2. If we multiply them we get an expression with all real coefficients, because it is self-conjugate. In addition, replacing $y$ with $wy$ doesn't change the product, so the powers of $y$ present in the result have to be divisible by 3. Analogously, the powers of $z$ have to be divisible by 3. Consequently, as the product is homogeneous, the powers of $x$ have to be divisible by 3.

Thus the resulting product is a polynomial of $a$, $b$, and $c$. Explicitly, if we multiply the denominator and numerator by $(x + w y + w^2 z) (x + w^2 y + w z) (x + y + w z) (x + y + w^2 z) (x + w  y + z) (x + w^2 y + z) (x + w y + w z) (x + w^2 y + w^2 z)$, the denominator becomes $(a+b + c)^3 - 27abc$.

We might try to simplify the calculation by initially multiplying the numerator and denominator by $(x + w y + w^2 z) (x + w^2 y + w z) = x^2 + y^2 - y z + z^2 - x (y + z)$. We get

$$\frac{1}{\sqrt[3]{a}+\sqrt[3]{b}+\sqrt[3]{c}} = \frac{\sqrt[3]{a}^2 + \sqrt[3]{b}^2 + \sqrt[3]{c}^2 - \sqrt[3]{ab} - \sqrt[3]{ac} - \sqrt[3]{bc}}{a + b + c - 3 \sqrt[3]{abc}}.$$

Now the denominator is of the form $s - \sqrt[3]{t}$ and has only one radical. To get rid of it we can multiply the denominator and the numerator by $s^2 + s\sqrt[3]{t} + \sqrt[3]{t}^2$. The denominator becomes $s^3 - t = (a+b+c)^3 - 27abc$.

\subsection*{Solution to problem~\ref{pr:68}}
\begin{figure}[htbp]
\begin{center}
\ifpdf
  \input{pr68.pdftex_t}
\else
  \input{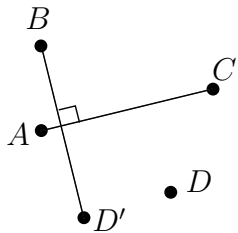}
\fi
\end{center}
\caption{Problem~\ref{pr:68}: Constructing a square from four given points.}
\label{pr68:fig}
\end{figure}

Suppose $A$, $B$, $C$ and $D$ denote the points on the sides of the square, in order. Connect $A$ and $C$. Draw a perpendicular from $B$ to the segment $AC$. Let point $D'$ lie on the perpendicular given that the distance $BD'$ is the same as the length of $AC$. Then $D'$ belongs to the same side (or extension of the side) of the square as $D$. See Figure~\ref{pr68:fig}.  Once we know from $DD'$ how the square is oriented, the rest is easy.

\subsection*{Solution to problem~\ref{pr:71}}
Let $A$ be the point on the function that lies on the middle line between and parallel to the two lines. If we move the cut off point to the right, then the area under the function increases faster than the area over the function decreases. Moving to the left is symmetric.

\end{document}

%% file: pr48.pdftex_t
\begin{picture}(0,0)%
\includegraphics{pr48.pdftex}%
\end{picture}%
\setlength{\unitlength}{3674sp}%
\begingroup\makeatletter\ifx\SetFigFont\undefined%
\gdef\SetFigFont#1#2{%
  \fontsize{#1}{#2pt}%
  \selectfont}%
\fi\endgroup%
\begin{picture}(2430,2578)(2761,-6941)
\put(5101,-5836){\makebox(0,0)[lb]{\smash{{\SetFigFont{12}{14.4}{\color[rgb]{0,0,0}$F$}%
}}}}
\put(3151,-5161){\makebox(0,0)[lb]{\smash{{\SetFigFont{12}{14.4}{\color[rgb]{0,0,0}$E$}%
}}}}
\put(3676,-4486){\makebox(0,0)[lb]{\smash{{\SetFigFont{12}{14.4}{\color[rgb]{0,0,0}$B$}%
}}}}
\put(5176,-4786){\makebox(0,0)[lb]{\smash{{\SetFigFont{12}{14.4}{\color[rgb]{0,0,0}$C$}%
}}}}
\put(4951,-6886){\makebox(0,0)[lb]{\smash{{\SetFigFont{12}{14.4}{\color[rgb]{0,0,0}$D$}%
}}}}
\put(2776,-5836){\makebox(0,0)[lb]{\smash{{\SetFigFont{12}{14.4}{\color[rgb]{0,0,0}$A$}%
}}}}
\end{picture}%

%% file: pr71.pdftex_t
\begin{picture}(0,0)%
\includegraphics{pr71.pdftex}%
\end{picture}%
\setlength{\unitlength}{6163sp}%
\begingroup\makeatletter\ifx\SetFigFont\undefined%
\gdef\SetFigFont#1#2{%
  \fontsize{#1}{#2pt}%
  \selectfont}%
\fi\endgroup%
\begin{picture}(953,873)(1792,-3164)
\end{picture}%

%% file: pr22.pdftex_t
\begin{picture}(0,0)%
\includegraphics{pr22.pdftex}%
\end{picture}%
\setlength{\unitlength}{2557sp}%
\begingroup\makeatletter\ifx\SetFigFont\undefined%
\gdef\SetFigFont#1#2{%
  \fontsize{#1}{#2pt}%
  \selectfont}%
\fi\endgroup%
\begin{picture}(4351,5557)(2963,-5281)
\put(6451,104){\makebox(0,0)[lb]{\smash{{\SetFigFont{25}{30.0}{\color[rgb]{0,0,0}$M$}%
}}}}
\end{picture}%

%% file: pr26.pdftex_t
\begin{picture}(0,0)%
\includegraphics{pr26.pdftex}%
\end{picture}%
\setlength{\unitlength}{3570sp}%
\begingroup\makeatletter\ifx\SetFigFont\undefined%
\gdef\SetFigFont#1#2{%
  \fontsize{#1}{#2pt}%
  \selectfont}%
\fi\endgroup%
\begin{picture}(2745,4643)(2161,-5106)
\put(3471,-3231){\makebox(0,0)[lb]{\smash{{\SetFigFont{12}{14.4}{\color[rgb]{0,0,0}$x$}%
}}}}
\put(3306,-2051){\makebox(0,0)[lb]{\smash{{\SetFigFont{12}{14.4}{\color[rgb]{0,0,0}$x$}%
}}}}
\put(3046,-3331){\makebox(0,0)[lb]{\smash{{\SetFigFont{12}{14.4}{\color[rgb]{0,0,0}$y$}%
}}}}
\put(3256,-2511){\makebox(0,0)[lb]{\smash{{\SetFigFont{12}{14.4}{\color[rgb]{0,0,0}$y$}%
}}}}
\put(3496,-586){\makebox(0,0)[lb]{\smash{{\SetFigFont{12}{14.4}{\color[rgb]{0,0,0}$B$}%
}}}}
\put(4891,-2851){\makebox(0,0)[lb]{\smash{{\SetFigFont{12}{14.4}{\color[rgb]{0,0,0}$C = B'$}%
}}}}
\put(2176,-2866){\makebox(0,0)[lb]{\smash{{\SetFigFont{12}{14.4}{\color[rgb]{0,0,0}$A$}%
}}}}
\put(3346,-3001){\makebox(0,0)[lb]{\smash{{\SetFigFont{12}{14.4}{\color[rgb]{0,0,0}$O'$}%
}}}}
\put(2921,-2131){\makebox(0,0)[lb]{\smash{{\SetFigFont{12}{14.4}{\color[rgb]{0,0,0}$O$}%
}}}}
\put(3496,-5051){\makebox(0,0)[lb]{\smash{{\SetFigFont{12}{14.4}{\color[rgb]{0,0,0}$C'$}%
}}}}
\end{picture}%

%% file: pr48sol1.pdftex_t
\begin{picture}(0,0)%
\includegraphics{pr48sol1.pdftex}%
\end{picture}%
\setlength{\unitlength}{5910sp}%
\begingroup\makeatletter\ifx\SetFigFont\undefined%
\gdef\SetFigFont#1#2{%
  \fontsize{#1}{#2pt}%
  \selectfont}%
\fi\endgroup%
\begin{picture}(2805,2578)(2611,-6941)
\put(4576,-5311){\makebox(0,0)[lb]{\smash{{\SetFigFont{12}{14.4}{\color[rgb]{0,0,0}$C'$}%
}}}}
\put(5401,-6136){\makebox(0,0)[lb]{\smash{{\SetFigFont{12}{14.4}{\color[rgb]{0,0,0}$D'$}%
}}}}
\put(2626,-5161){\makebox(0,0)[lb]{\smash{{\SetFigFont{12}{14.4}{\color[rgb]{0,0,0}$A' = B' =$}%
}}}}
\put(5101,-5836){\makebox(0,0)[lb]{\smash{{\SetFigFont{12}{14.4}{\color[rgb]{0,0,0}$F$}%
}}}}
\put(3226,-5161){\makebox(0,0)[lb]{\smash{{\SetFigFont{12}{14.4}{\color[rgb]{0,0,0}$E$}%
}}}}
\put(3676,-4486){\makebox(0,0)[lb]{\smash{{\SetFigFont{12}{14.4}{\color[rgb]{0,0,0}$B$}%
}}}}
\put(5176,-4786){\makebox(0,0)[lb]{\smash{{\SetFigFont{12}{14.4}{\color[rgb]{0,0,0}$C$}%
}}}}
\put(2851,-5836){\makebox(0,0)[lb]{\smash{{\SetFigFont{12}{14.4}{\color[rgb]{0,0,0}$A$}%
}}}}
\put(4951,-6886){\makebox(0,0)[lb]{\smash{{\SetFigFont{12}{14.4}{\color[rgb]{0,0,0}$D$}%
}}}}
\end{picture}%

%% file: pr10.pdftex_t
\begin{picture}(0,0)%
\includegraphics{pr10.pdftex}%
\end{picture}%
\setlength{\unitlength}{3167sp}%
\begingroup\makeatletter\ifx\SetFigFont\undefined%
\gdef\SetFigFont#1#2{%
  \fontsize{#1}{#2pt}%
  \selectfont}%
\fi\endgroup%
\begin{picture}(5235,3583)(2063,-3716)
\put(3488,-256){\makebox(0,0)[lb]{\smash{{\SetFigFont{12}{14.4}{\color[rgb]{0,0,0}$B$}%
}}}}
\put(5933,-2153){\makebox(0,0)[lb]{\smash{{\SetFigFont{12}{14.4}{\color[rgb]{0,0,0}$B'$}%
}}}}
\put(6676,-2791){\makebox(0,0)[lb]{\smash{{\SetFigFont{12}{14.4}{\color[rgb]{0,0,0}$M'$}%
}}}}
\put(7283,-3391){\makebox(0,0)[lb]{\smash{{\SetFigFont{12}{14.4}{\color[rgb]{0,0,0}$C$}%
}}}}
\put(2078,-3503){\makebox(0,0)[lb]{\smash{{\SetFigFont{12}{14.4}{\color[rgb]{0,0,0}$A$}%
}}}}
\put(5476,-2543){\makebox(0,0)[lb]{\smash{{\SetFigFont{12}{14.4}{\color[rgb]{0,0,0}$K'$}%
}}}}
\put(5266,-3661){\makebox(0,0)[lb]{\smash{{\SetFigFont{12}{14.4}{\color[rgb]{0,0,0}$A'$}%
}}}}
\put(2926,-1261){\makebox(0,0)[lb]{\smash{{\SetFigFont{12}{14.4}{\color[rgb]{0,0,0}$K$}%
}}}}
\put(2401,-2611){\makebox(0,0)[lb]{\smash{{\SetFigFont{12}{14.4}{\color[rgb]{0,0,0}$L$}%
}}}}
\end{picture}%

%% file: pr19.pdftex_t
\begin{picture}(0,0)%
\includegraphics{pr19.pdftex}%
\end{picture}%
\setlength{\unitlength}{4778sp}%
\begingroup\makeatletter\ifx\SetFigFont\undefined%
\gdef\SetFigFont#1#2{%
  \fontsize{#1}{#2pt}%
  \selectfont}%
\fi\endgroup%
\begin{picture}(2974,2719)(4711,-4300)
\put(6076,-2236){\makebox(0,0)[lb]{\smash{{\SetFigFont{12}{14.4}{\color[rgb]{0,0,0}$A$}%
}}}}
\put(6076,-3811){\makebox(0,0)[lb]{\smash{{\SetFigFont{12}{14.4}{\color[rgb]{0,0,0}$B$}%
}}}}
\put(4726,-3061){\makebox(0,0)[lb]{\smash{{\SetFigFont{12}{14.4}{\color[rgb]{0,0,0}$C$}%
}}}}
\put(5251,-1725){\makebox(0,0)[lb]{\smash{{\SetFigFont{12}{14.4}{\color[rgb]{0,0,0}$M$}%
}}}}
\end{picture}%

%% file: pr30i.pdftex_t
\begin{picture}(0,0)%
\includegraphics{pr30i.pdftex}%
\end{picture}%
\setlength{\unitlength}{6446sp}%
\begingroup\makeatletter\ifx\SetFigFont\undefined%
\gdef\SetFigFont#1#2{%
  \fontsize{#1}{#2pt}%
  \selectfont}%
\fi\endgroup%
\begin{picture}(1482,1837)(2209,-5375)
\put(3001,-3661){\makebox(0,0)[lb]{\smash{{\SetFigFont{12}{14.4}{\color[rgb]{0,0,0}$B$}%
}}}}
\put(3676,-5161){\makebox(0,0)[lb]{\smash{{\SetFigFont{12}{14.4}{\color[rgb]{0,0,0}$C$}%
}}}}
\put(3001,-5320){\makebox(0,0)[lb]{\smash{{\SetFigFont{12}{14.4}{\color[rgb]{0,0,0}$D$}%
}}}}
\put(2224,-5161){\makebox(0,0)[lb]{\smash{{\SetFigFont{12}{14.4}{\color[rgb]{0,0,0}$A$}%
}}}}
\end{picture}%

%% file: pr30.pdftex_t
\begin{picture}(0,0)%
\includegraphics{pr30.pdftex}%
\end{picture}%
\setlength{\unitlength}{2248sp}%
\begingroup\makeatletter\ifx\SetFigFont\undefined%
\gdef\SetFigFont#1#2{%
  \fontsize{#1}{#2pt}%
  \selectfont}%
\fi\endgroup%
\begin{picture}(6322,2542)(1939,-5831)
\end{picture}%

%% file: pr42a.pdftex_t
\begin{picture}(0,0)%
\includegraphics{pr42a.pdftex}%
\end{picture}%
\setlength{\unitlength}{1464sp}%
\begingroup\makeatletter\ifx\SetFigFont\undefined%
\gdef\SetFigFont#1#2{%
  \fontsize{#1}{#2pt}%
  \selectfont}%
\fi\endgroup%
\begin{picture}(7830,8086)(286,-6950)
\put(3001,-6886){\makebox(0,0)[rb]{\smash{{\SetFigFont{12}{14.4}{\color[rgb]{0,0,0}$B: (-7, -24)$}%
}}}}
\put(301,-3286){\makebox(0,0)[rb]{\smash{{\SetFigFont{12}{14.4}{\color[rgb]{0,0,0}$A: (-25,0)$}%
}}}}
\put(8101,-3286){\makebox(0,0)[lb]{\smash{{\SetFigFont{12}{14.4}{\color[rgb]{0,0,0}$C: (25, 0)$}%
}}}}
\put(3001,989){\makebox(0,0)[rb]{\smash{{\SetFigFont{12}{14.4}{\color[rgb]{0,0,0}$(-7, 24)$}%
}}}}
\put(5401,989){\makebox(0,0)[lb]{\smash{{\SetFigFont{12}{14.4}{\color[rgb]{0,0,0}$(7,24)$}%
}}}}
\put(5401,-6886){\makebox(0,0)[lb]{\smash{{\SetFigFont{12}{14.4}{\color[rgb]{0,0,0}$(-7,24)$}%
}}}}
\put(3001,-2611){\makebox(0,0)[rb]{\smash{{\SetFigFont{12}{14.4}{\color[rgb]{0,0,0}$D$}%
}}}}
\end{picture}%

%% file: pr48sol2.pdftex_t
\begin{picture}(0,0)%
\includegraphics{pr48sol2.pdftex}%
\end{picture}%
\setlength{\unitlength}{6209sp}%
\begingroup\makeatletter\ifx\SetFigFont\undefined%
\gdef\SetFigFont#1#2{%
  \fontsize{#1}{#2pt}%
  \selectfont}%
\fi\endgroup%
\begin{picture}(2670,2578)(2521,-6941)
\put(4351,-6136){\makebox(0,0)[lb]{\smash{{\SetFigFont{12}{14.4}{\color[rgb]{0,0,0}$C''$}%
}}}}
\put(2536,-5836){\makebox(0,0)[lb]{\smash{{\SetFigFont{12}{14.4}{\color[rgb]{0,0,0}$B'' =$}%
}}}}
\put(3676,-4486){\makebox(0,0)[lb]{\smash{{\SetFigFont{12}{14.4}{\color[rgb]{0,0,0}$B$}%
}}}}
\put(5176,-4786){\makebox(0,0)[lb]{\smash{{\SetFigFont{12}{14.4}{\color[rgb]{0,0,0}$C$}%
}}}}
\put(2851,-5836){\makebox(0,0)[lb]{\smash{{\SetFigFont{12}{14.4}{\color[rgb]{0,0,0}$A$}%
}}}}
\put(4951,-6886){\makebox(0,0)[lb]{\smash{{\SetFigFont{12}{14.4}{\color[rgb]{0,0,0}$D$}%
}}}}
\end{picture}%

%% file: pr49.pdftex_t
\begin{picture}(0,0)%
\includegraphics{pr49.pdftex}%
\end{picture}%
\setlength{\unitlength}{3450sp}%
\begingroup\makeatletter\ifx\SetFigFont\undefined%
\gdef\SetFigFont#1#2{%
  \fontsize{#1}{#2pt}%
  \selectfont}%
\fi\endgroup%
\begin{picture}(1986,1708)(2096,-3940)
\end{picture}%

%% file: pr50.pdftex_t
\begin{picture}(0,0)%
\includegraphics{pr50.pdftex}%
\end{picture}%
\setlength{\unitlength}{5107sp}%
\begingroup\makeatletter\ifx\SetFigFont\undefined%
\gdef\SetFigFont#1#2{%
  \fontsize{#1}{#2pt}%
  \selectfont}%
\fi\endgroup%
\begin{picture}(1150,950)(2213,-3251)
\end{picture}%

%% file: pr68.pdftex_t
\begin{picture}(0,0)%
\includegraphics{pr68.pdftex}%
\end{picture}%
\setlength{\unitlength}{4740sp}%
\begingroup\makeatletter\ifx\SetFigFont\undefined%
\gdef\SetFigFont#1#2{%
  \fontsize{#1}{#2pt}%
  \selectfont}%
\fi\endgroup%
\begin{picture}(1133,1239)(2056,-3337)
\put(3009,-3075){\makebox(0,0)[lb]{\smash{{\SetFigFont{12}{14.4}{\color[rgb]{0,0,0}$D$}%
}}}}
\put(3147,-2491){\makebox(0,0)[lb]{\smash{{\SetFigFont{12}{14.4}{\color[rgb]{0,0,0}$C$}%
}}}}
\put(2170,-2221){\makebox(0,0)[lb]{\smash{{\SetFigFont{12}{14.4}{\color[rgb]{0,0,0}$B$}%
}}}}
\put(2071,-2829){\makebox(0,0)[lb]{\smash{{\SetFigFont{12}{14.4}{\color[rgb]{0,0,0}$A$}%
}}}}
\put(2521,-3282){\makebox(0,0)[lb]{\smash{{\SetFigFont{12}{14.4}{\color[rgb]{0,0,0}$D'$}%
}}}}
\end{picture}%